\documentclass[twocolumn]{autart}    

\usepackage{cite}
\usepackage{placeins}
\usepackage{amsmath,amssymb}
\usepackage{algorithmic}
\usepackage{graphicx}
\usepackage{textcomp}
\usepackage{mathtools}
\usepackage{enumerate}
\usepackage{graphicx,xcolor,cases}          
\newtheorem{theorem}{Theorem}
\newtheorem{lemma}[theorem]{Lemma}                               
\begin{document}

\begin{frontmatter}

\title{Finite-time stabilization of ladder multi-level quantum systems} 
\thanks[footnoteinfo]{
This work was supported in part by the National Natural Science Foundation of China under Grant 62373342. 
}
\author[AD1]{Zeping Su}\ead{zeping.su@mail.ustc.edu.cn}, 
\author[AD1]{Sen Kuang\corauthref{correspondingauthor}}
\corauth[correspondingauthor]{Corresponding author.}
\ead{skuang@ustc.edu.cn},
\author[AD2]{Daoyi Dong}\ead{daoyidong@gmail.com}                   
\address[AD1]{Department of Automation, University of Science and Technology of China, Hefei 230027, PR China}  
\address[AD2]{Australian Artificial Intelligence Institute, Faculty of Engineering and IT, University of Technology of Sydney, NSW 2007, Australia} 



\begin{abstract}        
In this paper, a novel continuous non-smooth control strategy is proposed to achieve finite-time stabilization of ladder quantum systems. We first design a universal fractional-order control law for a ladder $n$-level quantum system using a distance-based Lyapunov function, and then apply the Filippov solution in the sense of differential inclusions and the LaSalle's invariance principle to prove the existence and uniqueness of the solution of the ladder system under the continuous non-smooth control law. Both asymptotic stability and finite-time stability for the ladder system is rigorously established by applying Lyapunov stability theory and finite-time stability criteria. We also derive an upper bound of the time required for convergence to an eigenstate of the intrinsic Hamiltonian. Numerical simulations on a rubidium ladder three-level atomic system validate the effectiveness of the proposed method.
\end{abstract}
\end{frontmatter}
\section{Introduction}
Quantum information technology is developing rapidly. The special properties of quantum systems, such as superposition and entanglement enable advanced information processing techniques, including quantum communication \cite{bennett1993teleporting} and quantum sensing \cite{bongs2023quantum}, which have the potential to surpass classical information technology \cite{awschalom2013quantum,albertini2003notions}. Quantum control is the foundation of quantum information technology. Many classical control methods are applied to quantum systems, including optimal control \cite{dolde2014high,PhysRevA.95.063418}, Lyapunov control \cite{KUANG2017164,GRIGORIU20122229}, sliding mode control \cite{Dong_2009}, $H^\infty$ control \cite{XIANG20178}, fault-tolerant control and filtering \cite{GAO2016125}, and learning and robust control \cite{8759071,9792196,DONG2023}, etc.

One crucial objective of quantum control is to achieve high control precision with a fast convergence rate under the action of control fields. Finite-time control offers superior control accuracy, rapid convergence, and strong robustness against various uncertainties \cite{doi:10.1137/S0363012997321358}. Kuang et al. \cite{KUANG2021109327} proposed a continuous non-smooth control law and achieved finite-time control for a spin-1/2 system. Building upon this work, Li et al.  \cite{LI2022100005} used fixed-time Lyapunov control technique and bi-limit homogeneity theory to stabilize a two-level quantum system to one eigenstate of its intrinsic Hamiltonian within a fixed time, regardless of the initial state. Based on relative state distance, Taslima et al. \cite{taslima2024lyapunov} designed a switching control law by considering prescribed-time Lyapunov control technique and stabilized a two-level quantum system to one eigenstate of its internal Hamiltonian within a specified time. Since this approach allows for the flexible selection of the upper bound of the stabilization time, the system is stabilized within the designated time without being affected by initial conditions. Wu et al. \cite{wu2024arbitrary} extended the controlled model to open two-level qubit systems with decoherence. To prevent the qubit from entering invariant or singular sets, they proposed a switching control that drives the system to a sufficiently small neighborhood of the target state, and meanwhile guarantees finite-time stability and finite-time contraction stability. The aforementioned methods are all restricted on two-level systems and cannot be directly extended to multi-level systems.

However, in many practical applications, it is necessary to consider multi-level quantum systems, e.g., in the scenarios of electromagnetically induced transparency (EIT) \cite{RevModPhys.77.633} and single-photon generation \cite{PhysRevLett.89.067901}.
Besides, multi-level systems are also widely applied in quantum computing, particularly in the preparation and measurement of qubits \cite{timoney2011quantum,PhysRevLett.111.140501}. In quantum sensing, quantum entanglement and coherence in high-dimensional systems significantly enhance measurement sensitivity and resolution \cite{huang2024entanglement}. Although one can consider optimal control for high-dimensional systems \cite{KHANEJA2005296}, it requires complex numerical calculations. In fact, it remains a significant challenge for achieving precise convergence of a high-dimensional system to the target state within finite time.

In this paper, we explore the finite-time control of multi-level ladder quantum systems. The ladder structure has potential for wide applications in quantum computing, quantum networks, and EIT\cite{PhysRevResearch.7.L012049}. A ladder quantum system has energy transitions that occur in a stepwise manner, where each level only connects to its adjacent one. This structure allows for precise control of quantum states and enables effective photon-state transfers, crucial for quantum computing and communication \cite{PhysRevA.75.013809}. By manipulating the lifetime ratio, the quantum interference of cascaded photons can be enhanced, which aids in quantum repeater development and quantum network advancement \cite{PhysRevLett.125.233605}. Additionally, the ladder quantum system using Rydberg atoms offers high-precision electric field detection and efficient state control through EIT, playing a key role in quantum sensing technologies \cite{PhysRevA.104.043103}. 

The contributions of this paper are as follows. First, for arbitrarily dimensional ladder quantum systems, we design a universal fractional-order control law and achieve convergenve to the target state within finite time. This approach extends the existing finite-time stabilization methods to ladder $n$-level quantum systems. Second, since high-dimensional systems result in a significant increase in the number of non-Lipschitz regions under the continuous non-smooth control law, we combine the Filippov solution based on differential inclusions with LaSalle's invariance principle to prove the existence and uniqueness of the dynamical solution of the system under the designed controller.
Third, we mathematically prove that the system under the designed control law is not only asymptotically stable but also finite-time convergent. Lastly, we derive an upper bound of the time required for the system to reach the target state.


\section{System Model and Control Problem}
The quantum state of an $n$-level system is represented by a unit vector $\left| \psi \right\rangle$ in the Hilbert space ${\mathbb{C}}^{n}$, and its evolution is governed by the Schrödinger equation
\begin{equation} \label{b1}
i\hbar|\dot{\psi}(t)\rangle = \Big( H_0 + \sum_{j=1}^{r} H_j u_j \Big) |\psi(t)\rangle,
\end{equation}
where $\hbar$ is the reduced Planck constant, which is set to 1 in this paper; ${u_j}$ is the external control field to be designed; ${H_0}$ is the intrinsic Hamiltonian of the system; and ${H_j}$ represents the control Hamiltonian induced by the interaction between ${u_j}$ and the system.

We denote the coherence vector of \( |\psi\rangle \) as $[\xi_1, \ldots, \xi_{n^2-1}]^T$ $=  [\langle \psi | \sigma_1 | \psi \rangle, \ldots, \langle \psi | \sigma_{n^2-1} | \psi \rangle]^T\triangleq s =\in \mathbb{R}^{n^2-1}$, where \( \{ \sigma_k \}_{k=0}^{n^2-1} \) forms an orthonormal basis of the space of \( n \times n \) complex Hermitian matrices with \( \sigma_0 = \frac{I_n}{\sqrt{n}} \). The set of all coherence vectors forms the Bloch space \( \mathcal{B}(\mathbb{R}^{n^2-1}) \) \cite{KIMURA2003339}. If the control \( u_k \) in \eqref{b1} is a continuous function of the state \( s \), then the system \eqref{b1} can be rewritten as\cite{KUANG2021109327}
\begin{equation} \label{ss}
\dot{s}(t) = f(s(t)), \quad s(t) \in \mathcal{B}(\mathbb{R}^{n^2-1})
\end{equation}
where \( f : \mathcal{B}(\mathbb{R}^{n^2-1}) \rightarrow \mathcal{B}(\mathbb{R}^{n^2-1}) \) is a continuous function defined on \( \mathcal{B}(\mathbb{R}^{n^2-1}) \). 

A ladder quantum system consists of energy levels arranged in sequence, where each energy level is coupled to its adjacent upper and lower energy levels. We use $n-1$ control fields to regulate the evolution of the system and assume that its intrinsic and control Hamiltonians have the following form
\begin{gather}
H_0 = \operatorname{diag}(\lambda_1, \lambda_2, \dots, \lambda_n),\label{b2}\\
H_p = i X_p, \quad (p=1, \ldots, n-1), \label{b3}
\end{gather}
where \( X_p \) has non-zero elements \(-1\) and \( 1 \) at the \( (p, p+1) \) position and \( (p+1, p) \) position, respectively, and all other elements are zero. 

In polar coordinates, we can express the state \( \left| \psi \right\rangle \) as
\begin{equation}\label{b4}
\left| \psi  \right\rangle =[{{r}_{1}}{{e}^{i{{\phi }_{1}}}},{{r}_{2}}{{e}^{i{{\phi }_{2}}}},\cdots ,{{r}_{n}}{{e}^{i{{\phi }_{n}}}}]
\end{equation}
with $r_1, \ldots, r_n \in [0,1]$ and $r_1^2 + \dots + r_n^2 = 1$. In particular, we define $\phi_j = 0$ for $r_j = 0$.
Thus, the system \eqref{b1} with the Hamiltonians \eqref{b2} and \eqref{b3} can be written as
\begin{equation}\label{b5}
\begin{cases}
\dot{r}_1 e^{i \phi_1} + i r_1 e^{i \phi_1} \dot{\phi}_1= -i \lambda_1 r_1 e^{i \phi_1} -  r_2 e^{i \phi_2} u_1 \\
\dot{r}_j e^{i \phi_j} + i r_j e^{i \phi_j} \dot{\phi}_j= -i \lambda_j r_j e^{i \phi_j} + r_{j-1} e^{i \phi_{j-1}} u_{j-1}\\
\hspace{3.2cm} -  r_{j+1} e^{i \phi_{j+1}} u_{j}, (j=2,\ldots ,n-1)\\
\dot{r}_n e^{i \phi_n}\! + i r_n e^{i \phi_n} \dot{\phi}_n
= -i \lambda_n r_n e^{i \phi_n}\! + r_{n-1} e^{i \phi_{n-1}} u_{n-1}.
\end{cases}
\end{equation}
Define the phase difference between ${{e}^{i{{\phi }_{i}}}}$ and ${{e}^{i{{\phi }_{j}}}}$ as ${{\phi }_{ij}}\triangleq {{\phi }_{i}}-{{\phi }_{j}}$. Then, \eqref{b5} is equivalent to
\begin{equation}\label{b6}
\begin{cases}
\dot{r}_1 = -u_1 r_2 \cos \phi_{21} \\
\dot{r}_j = u_{j-1} r_{j-1} \cos \phi_{j(j-1)}- u_j r_{j+1} \cos \phi_{(j+1)j} \\
\dot{r}_n = u_{n-1} r_{n-1} \cos \phi_{n(n-1)}
\end{cases}
\end{equation}
and
\begin{equation}\label{b7}
\begin{cases}
r_1 \dot{\phi}_1 = -\lambda_1 r_1 - u_1 r_2 \sin \phi_{21} \\
r_j \dot{\phi}_j = -\lambda_j r_j - u_{j-1} r_{j-1} \sin \phi_{j(j-1)} \\
\qquad\quad\, - u_j r_{j+1} \sin \phi_{(j+1)j} \\
r_n \dot{\phi}_n = -\lambda_n r_n - u_{n-1} r_{n-1} \sin \phi_{n(n-1)}.
\end{cases}
\end{equation}
The control task in this paper is to design the control laws in \eqref{b6} and \eqref{b7} to achieve exact convergence to an eigenstate \( \left| \psi_f \right\rangle \) of \( H_0 \) in finite time. Without loss of generality, we assume that the target state is the last eigenstate of \( H_0 \), i.e., \( | \psi_f \rangle = | n-1 \rangle = [0, 0, \dots, 0, 1]^T \).

\section{Finite-Time Controller Design}
Let us consider the following Lyapunov function \cite{CONG200728}:
\begin{equation} \label{e2}
V = 1 - \left| \langle \psi_f | \psi \rangle \right|^2.
\end{equation}
Calculating the first-order time derivative of $V$ yields
\begin{equation}\label{e3}
\dot{V}= - 2 \sum_{j=1}^{n-1} u_j \left| \langle\psi | \psi_f \rangle \right| \Im \Big[e^{i \angle \langle \psi | \psi_f \rangle} \langle \psi_f | H_j | \psi \rangle \Big],
\end{equation}
where $\angle\langle \psi | \psi_f \rangle$ denotes the argument of $\langle \psi | \psi_f \rangle$. For $\langle \psi | \psi_f \rangle = 0$, we define $\angle \langle \psi | \psi_f \rangle = 0^\circ$.

By using \eqref{b4}, \eqref{e2} can be written as
\begin{equation} \label{e6}
V = 1 - \left| \langle \psi_f | \psi \rangle \right|^2 = r_{1}^{2}+r_{2}^{2}+\cdots +r_{n-1}^{2}.
\end{equation}
With \eqref{b2} and \eqref{b3}, $\langle \psi_f | H_k | \psi \rangle = 0$ holds for $k=1,\ldots,n-1$. Thus, $\dot{V}$ in \eqref{e3} can be calculated as
\begin{equation}\label{e7}
\dot{V}=-2{{u}_{n-1}}{{r}_{n}}{{r}_{n-1}}\cos {{\phi }_{n(n-1)}}.
\end{equation}
To ensure that $\dot{V} \leq 0$, we design ${{u}_{j}}$ as the following continuous non-smooth form:
\begin{equation}\label{e8}
u_{j} = k_{j} \operatorname{sign}(r_{j} \cos \phi_{(j+1)j}) \left| r_{j} \cos \phi_{(j+1)j} \right|^{\alpha_{j}}
\end{equation}
with $k_{j} > 0$, ${{\alpha }_{j}} \in (0,1)$, and $j=1,\ldots ,n-1$. 

Substituting \eqref{e8} into \eqref{e7}, we obtain
\begin{equation}\label{e9}
\dot{V}=-2{{k}_{n-1}}{{r}_{n}}|{{r}_{n-1}}\cos {{\phi }_{n(n-1)}}{{|}^{{{\alpha }_{n-1}}+1}}.
\end{equation}
%
It should be noted that the dynamic system in \eqref{b6} and \eqref{b7} with the controller \eqref{e8} is not globally Lipschitz continuous. Here, we provide a sufficient condition to show the existence and uniqueness of solutions to the system. 
\begin{thm}\label{Theorem1}
With the controller \eqref{e8}, the non-degenerate $n$-level quantum system in \eqref{b6} and \eqref{b7}, equivalently, the system \eqref{b1} with the Hamiltonians \eqref{b2} and \eqref{b3}, has a unique continuously differentiable solution, i.e., the equivalence class of the target state, for any initial state.
\end{thm}
\textbf{Proof.} Since the control law \( u_{j} \) is continuous but non-smooth, the system in \eqref{b6} and \eqref{b7} is not globally Lipschitz continuous, which may lead to the cases of no solution and non-unique solutions. Nevertheless, it still has a unique solution under the controller \eqref{e8}. The proof is carried out in the following two steps.
\begin{enumerate}[Step 1]
  \item The differential inclusion system \eqref{e16} below has a Filippov solution for each initial state, and each Filippov solution corresponds to a solution of the system \eqref{b1}. Thus, the system \eqref{b1} has a solution.
  \item Using LaSalle's invariance principle, we prove the uniqueness of the solution to the system \eqref{b1}.
\end{enumerate}

\textbf{Proof of Step 1. }For the quantum system \eqref{b1}, we construct a set-valued mapping \( \mathcal{F} \) in the sense of the Filippov extension:
\begin{equation}\label{e16}
\dot{\left| \psi  \right\rangle}\in \mathcal{F}(\left| \psi  \right\rangle):=\underset{\delta >0}{\mathop{\bigcap }}\,\mathrm{clco}f(\left| \psi  \right\rangle+\delta \mathcal{B}),
\end{equation}
where $\mathrm{clco}(\mathcal{A})$ denotes the closed convex hull of $\mathcal{A}$; $\mathcal{B}$ is the unit ball in ${{\mathbb{R}}^{d}}$; and $f$ corresponds to the Schrödinger equation governing the system dynamics, i.e.,
\[
f(|\psi (t)\rangle )=|\dot{\psi}(t)\rangle=-i\Big( {{H}_{0}}+\sum\limits_{j=1}^{n-1}{{{H}_{j}}}{{u}_{j}} \Big)|\psi (t)\rangle. 
\]
The definition of $\mathcal{F}(\left| \psi \right\rangle)$ implies that it is a non-empty, bounded, closed, and convex set\cite{GRIGORIU20122229}. Since $H_0$, $H_1$ and $u_k$ are bounded, $\dot{| {{{{\psi }}}_{1}} \rangle} = f(| {{\psi }_{1}}\rangle )$ is also bounded. Although \( f \) is continuous non-smooth (due to the continuous non-smooth nature of \( u_k \)), the Schrödinger equation itself remains a single-valued function. Thus, based on Lemma~\ref{Lemma2} in the Appendix, the following equation has a solution, which corresponds to the Filippov extension of \eqref{e16}:
\[
|\dot{\psi}\rangle = \mathcal{F}(|\psi\rangle), \quad |\psi(0)\rangle = |\psi_0\rangle \in \mathbb{R}^n.
\]
Furthermore, since \( f \) is continuous everywhere, $\left| {{\psi }_{1}} \right\rangle $ is a solution of \( \mathcal{F} \) and also a solution of \( f \). That is, for any initial state $\left| {{\psi }_{0}} \right\rangle $, every Filippov solution $\left| {{\psi }_{1}} \right\rangle $ of the differential inclusion system \eqref{e16} corresponds to a solution of the system \eqref{b1}. Thus, the system \eqref{b1} has a solution, which may be not unique.

\textbf{Proof of Step 2. }According to \eqref{e9}, $\dot{V} = 0$ implies ${r}_{n} = 0$ or ${u}_{n-1} = 0$.

\textbf{Case 1}:  \( r_{n} = 0 \). In this case, \( \langle \psi | \psi_f \rangle = 0 \) and \( \phi_{n} = 0 \) hold. If \( u_{n-1} \neq 0 \), then \eqref{b6} implies \( \dot{r}_{n} \neq 0 \), meaning that the system is transferring toward the target state \( |\psi_{f} \rangle \). If \( u_{n-1} = 0 \), then it is included in one of the following Cases.

\textbf{Case 2}: ${{u}_{n-1}}=0$ and ${{u}_{n-2}}\ne 0$. In this case, we can consider several subcases.
\begin{itemize}
\item[(i)]When ${r}_{n-1} = 0$ and $\cos {{\phi }_{n(n-1)}} \ne 0$, \eqref{b6} implies ${{\dot{r}}_{n-1}} \ne 0$, and so ${r}_{n-1}$ will evolve to a non-zero value, which results in ${{u}_{n-1}} \ne 0$. It follows from \eqref{b6} that ${{\dot{r}}_{n}} > 0$ holds under the action of the controller \eqref{e8}, i.e., the system is moving toward the equivalence class of $\left| {{\psi }_{f}} \right\rangle$.
\item[(ii)]When ${r}_{n-1} \ne 0$, $\cos {{\phi }_{n(n-1)}} = 0$, and ${r}_{n} \ne 0$, the time derivative of $\cos \phi_{n(n-1)}$ is given by
\begin{equation}\label{x1}
( \cos \phi_{n(n-1)})' = -\sin \phi_{n(n-1)} ( \dot{\phi}_n - \dot{\phi}_{n-1}).
\end{equation}
Thus, according to \eqref{b7}, we have
\begin{numcases}{}
\dot{\phi}_{n-1} \!=\! -\lambda_{n-1} - u_{n-2} r_{n-2} \sin \phi_{(n-1)(n-2)} / r_{n-1}\label{e10}\\
\dot{\phi}_{n} = -\lambda_{n}.\label{e11}
\end{numcases}
Since $\cos \phi_{n(n-1)} = 0$, $\sin \phi_{n(n-1)} \neq 0$ holds. From \eqref{e10} and \eqref{e11}, it can be observed that $\dot{\phi}_n - \dot{\phi}_{n-1} \neq 0$ for most of the time. If $\dot{\phi}_n - \dot{\phi}_{n-1} = 0$, then $\ddot{\phi}_n = 0$. Since \( u_{n-2} \neq 0 \) and the relative phase keeps evolving, \( \ddot{\phi}_{n-1} \) varies continuously, which will make the system evolve to \( \dot{\phi}_n - \dot{\phi}_{n-1} \neq 0 \). 
In a word, \eqref{x1} implies that the system will evolve toward a state with $\cos \phi_{n(n-1)} \neq 0$ and so $u_{n-1} \neq 0$, while \eqref{e10} guarantees $\dot{r}_n > 0$. Therefore, the system is advancing toward the equivalence class of $| \psi_f\rangle$.

\item[(iii)]When ${r}_{n-1} \ne 0$, $\cos {{\phi }_{n(n-1)}} = 0$, and ${r}_{n} = 0$, \eqref{b7} means that
\begin{numcases}{}
\dot{\phi}_{n-1} \!=\! -\lambda_{n-1} - u_{n-2} r_{n-2} \sin \phi_{(n-1)(n-2)} / r_{n-1}\label{e12}\\
\phi_{n} = 0.\label{e13}
\end{numcases}
From \eqref{e13}, \( \dot{\phi}_n = 0 \) holds at all moments associated with \( r_n = 0 \). Thus, we have \( \cos {{\phi }_{n(n-1)}} = \cos {{\phi }_{n-1}} \). With \eqref{e12} and \( u_{n-2} \neq 0 \), and under the continuous evolution of the relative phase, we can obtain \( \dot{\phi}_{n-1} \neq 0 \). Equation \eqref{x1} implies \( \cos \dot{\phi}_{n-1} \neq 0 \), i.e., the system will evolve to \( \cos \phi_{n(n-1)} = \cos \phi_{n-1} \neq 0 \) and so \( u_{n-1}\ne 0 \). It is easily known from \eqref{e10} that \( \dot{r}_n > 0 \) and the system is progressing toward the equivalence class of \(| \psi_f\rangle \).
\item[(iv)]When \( r_{n-1} = 0 \), \( \cos \phi_{n(n-1)} = 0 \), and \( r_{n} \neq 0 \), it follows from \( u_{n-2} \neq 0 \) that both \( r_{n-2} \neq 0 \) and \( \cos \phi_{(n-1)(n-2)} \neq 0 \). Thus, according to \eqref{b6}, we have
\begin{equation}\label{x2}
\dot{r}_{n-1} = | r_{n-2} \cos \phi_{(n-1)(n-2)}|^{\alpha_{n-2}+1} \neq 0.
\end{equation}
This implies that \( r_{n-1} \) will evolve to a non-zero value, which is the subcase (ii) of Case 2. Ultimately, the system will evolve toward the equivalence class of \( |\psi_{f} \rangle \).
\item[(v)]When \( r_{n-1} = 0 \), \( \cos \phi_{n(n-1)} = 0 \), and \( r_{n} = 0 \), \eqref{x2} ensures ${{\dot{r}}_{n-1}} \ne 0$, i.e., \( r_{n-1} \) will evolve to a non-zero value, which is the subcase (iii) of Case 2, and so the system will evolve toward the equivalence class of \( |\psi_{f} \rangle \).
\end{itemize}
To sum up, for all subcases of Case 2, the system will transfer toward the equivalence class of \( |\psi_{f} \rangle \).

\textbf{Case $l$ ($l=3,4,\cdots n-1$)}:  ${{u}_{n-1}} = {{u}_{n-2}} = \cdots = {{u}_{n-l+1}} = 0$ and ${{u}_{n-l}} \ne 0$.
For this case, the detailed analysis is analogous to that of Case 2. All subcases of Case \( l \) will lead to \( {{u}_{n-l+1}} \neq 0 \), and so return to Case \( l-1 \). Thus, all cases will eventually reduce to Case 2, resulting in a transition toward the equivalence class of the target state. 

\textbf{Case $n$}: ${{u}_{n-1}} = {{u}_{n-2}} = \cdots = {{u}_{1}} = 0$.
\begin{itemize}
\item[(i)]When \( V = 0 \), the system reaches the target state \( |\psi_{f} \rangle \) and stabilises in its equivalence class. 
\item[(ii)]When $V \ne 0$, according to \eqref{b7}, we have ${{r}_{m}}{{\dot{\phi }}_{m}}=-{{\lambda }_{m}}{{r}_{m}},\,(m=1,\ldots,n)$. In fact, \( V \neq 0 \) means that there exists a set \( M \) such that \( r_{j} \neq 0 \), \(\exists j \in M \subseteq \{1, \ldots, n-1\} \). While \( r_{j} \neq 0 \), \( u_j = 0 \) can attribute to the phase condition \( \cos \phi_{(j+1)j} = 0 \). From \eqref{b7}, we can obtain
\begin{equation}\label{e22}
{{\dot{\phi }}_{j}}=-{{\lambda }_{j}},\ \exists j\in M.
\end{equation}
If \( j \in M \) and \( j+1 \in M \), i.e., \( r_{j} \neq 0 \) and \( r_{j+1} \neq 0 \), then it follows from \eqref{e22} that $\frac{d}{dt} ( \cos \phi_{(j+1)j}) \neq 0$ since \( H_{0} \) is non-degenerate. As a result, the relative phases will continue to evolve, leading to a non-zero control law \( u_{j} \) that drives the system back to the case analysed previously. If \( j \in M \) but \( j+1 \notin M \), i.e., \( r_{j} \neq 0 \) and \( r_{j+1} = 0 \), then we have \( \cos \phi_{(j+1)j} = \cos \phi_{j} \). Before \( r_{j+1} \neq 0 \) is satisfied, it follows from \eqref{e22} that
$
\frac{d}{dt}( \cos {{\phi }_{(j+1)j}}) = \frac{d}{dt}( \cos {{\phi }_{j}}) \neq 0
$.
Thus, the relative phases will evolve, resulting in a non-zero control \( u_{j} \), which drives the system back to one of the Cases analyzed previously. Then, the system will move toward the equivalence class of the target state \( |\psi_{f} \rangle \).
\end{itemize}
From the above, we know that when \( \dot{V} = 0 \), the system will move toward the equivalence class of the target state \( |\psi_{f} \rangle \) until it is reached. When \( \dot{V} \neq 0 \), the system will also move toward the equivalence class of the target state \( |\psi_{f} \rangle \) due to the effect of the controller. Therefore, the only solution of the system is the equivalence class of the target state \( |\psi_{f} \rangle \).

\section{Finite-Time Stability Analysis}
\subsection{Asymptotic Stability}
\begin{thm}\label{Theorem2}
Under the controller \eqref{e8}, the system \eqref{b1} with the Hamiltonians \eqref{b2} and \eqref{b3} is globally asymptotically stable, i.e., the system will be stabilised in the equivalence class of the target state \( |\psi_f\rangle = |n-1\rangle \).
\end{thm}
\textbf{Proof.}
According to Theorem 1 in \cite{CONG200728}, the largest invariant set of the system \eqref{b1} under the action of the controller \eqref{e8} is
$
{{S}^{2N-1}} \bigcap E
$
with $E =\{ |\psi\rangle : \langle \psi_f | H_k | \psi \rangle = 0,$ $ k = 1, \ldots, n-1\}$. For $ e^{i \phi_{n-1}}\ne 0$,  $\langle \psi_f | H_{n-1} | \psi \rangle = i r_{n-1} e^{i \phi_{n-1}} = 0$ means that
\begin{numcases}{}
r_{n-1} = 0,\label{ r_{n-1}=0}\\
u_{n-1} = 0.
\end{numcases}
Let us consider a function $m_{d(d-1)} = r_d r_{d-1} (d = 3, \ldots, n-1)$. Since $r_{d}$ is bounded, $m_{d(d-1)}$ is bounded too. Its first-order derivative is
\begin{equation}\label{dm}
\begin{aligned}
\dot{m}_{d(d-1)} &= \dot{r}_d r_{d-1} + r_d \dot{r}_{d-1}\\
&=\! u_{d-1} r_{d-1}^2 \cos \phi_{d(d-1)}\! -\! u_d r_{d+1} r_{d-1} \cos \phi_{(d+1)d}\\
&+\! u_{d-2} r_{d-2} r_d \cos \phi_{(d-1)(d-2)} \!-\! u_{d-1} r_d^2 \cos \phi_{d(d-1)}.
\end{aligned}
\end{equation}

\textbf{Step 1. } For \eqref{dm}, we let \( d = n-1 \) and define
\[
\begin{aligned}
g_{(n-1)1}(t) &= u_{n-2} r_{n-2}^2 \cos \phi_{(n-1)(n-2)} \\
&= k_{n-2} \left| \cos \phi_{(n-1)(n-2)} \right|^{1+\alpha_{n-2}} \times r_{n-2}^{2+\alpha_{n-2}},
\end{aligned}
\]
\[
\begin{aligned}
g_{(n-1)2}(t) &= -u_{n-1} r_n r_{n-2} \cos \phi_{n(n-1)} \\
&\quad + u_{n-3} r_{n-3} r_{n-1} \cos \phi_{(n-2)(n-3)} \\
&\quad - u_{n-2} r_{n-1}^2 \cos \phi_{(n-1)(n-2)}.
\end{aligned}
\]
With \eqref{ r_{n-1}=0}, we have \( \lim_{t \to \infty} g_{(n-1)2}(t) = 0 \). Since the first-order derivative of \( g_{(n-1)1} \) is bounded, \( g_{(n-1)1} \) is uniformly continuous. Thus, according to Lemma~\ref{Lemma3} in the Appendix, we can obtain \( \lim_{t \to \infty} g_{(n-1)1}(t) = 0 \). Hence, \( r_{n-2} = 0 \) or $\cos \phi_{(n-1)(n-2) =0}$. Based on Theorem~\ref{Theorem1}, 
if \( r_{n-2} \neq 0 \), then \( \cos \phi_{(n-1)(n-2)} \) cannot stabilize at 0. This means that \( r_{n-2} = 0 \).

\textbf{Step $k$.} Let $d = n-k$ in \eqref{dm}, where $k = 2,\ldots,n-2$. Similar to the proof of Step 1, by iterating backwards sequentially, we can obtain
\begin{equation}\label{r_n-j=0}
\begin{cases}
u_{n-k-1} = 0,\\ 
r_{n-k-1} = 0.
\end{cases}
\end{equation}
For the system \eqref{b1} with the Hamiltonians \eqref{b2} and \eqref{b3} and under the controller \eqref{e8}, \( r_k = 0 \), \( (k = 1, \ldots, n-1) \) holds. Thus, the largest invariant set of the system\eqref{b1} is the equivalence class \( e^{i\phi} |\psi_f\rangle \) of the target state \( |\psi_f\rangle \). According to Theorem 1 in \cite{CONG200728}, the control system is asymptotically stable \cite{CONG200728}. Since we do not impose any restrictions on initial states in the proof, asymptotic stability holds for any initial state. Therefore, the system \eqref{b1} is globally asymptotically stable in the sense of the equivalence class of \( |\psi_f\rangle \).

\subsection{Finite-time Stability}
\begin{lemma}\label{Lemma1}
Suppose that $\sum_{j=1}^{n}{r_{j}^{2}}=1$, where ${{r}_{j}}\ge 0$. Then, for $\forall \alpha \in (0,1)$, the following inequality holds
\end{lemma}
\begin{equation}\label{4lemma}
{{\Big( \sum\limits_{j=1}^{n}{r_{j}^{2}} \Big)}^{\frac{\alpha +1}{2}}}\le \sum\limits_{j=1}^{n}{r_{j}^{\alpha +1}}.
\end{equation}
\textbf{Proof.}
When ${{r}_{j}}=0$, $(j=1,\ldots,n)$, the inequality naturally holds. Otherwise, we denote $V_1 = \sum_{j=1}^{n} r_j^2$ and $V_2 = \sum_{j=1}^{n} r_j^{\alpha + 1}$. Since \( V_1 \) and \( V_2 \) are both positive definite functions with the homogeneity degrees of \( l_1 = 2 \) and \( l_2 = \alpha + 1 \), respectively, it follows from Lemma~\ref{Lemma4} in the Appendix that
\begin{equation}\label{4inequality}
\Big( \min_{\{Z: V_1(Z) = 1\}} V_2(Z) \Big) \big( V_1(r) \big)^{\frac{l_2}{l_1}} \leq V_2(r),
\end{equation}
where $z_w$ is the component of $Z$. Now, we minimize $V_2(Z)$, subject to the constraint $V_1(Z) = 1$. This optimization problem can be solved by using the Lagrange multiplier method. We construct a Lagrangian function as
\begin{equation}\label{e25}
\mathcal{L} = \sum_{w=1}^{n} z_w^{\alpha+1} - \lambda \Big( \sum_{w=1}^{n} z_w^2 - 1 \Big).
\end{equation}
It should be noted that when \eqref{e25} is used, the condition $z_w \geq 0$ always holds. Calculating the derivative with respect to $z_w$ and letting it be zero yields
\begin{equation}\label{e26}
\frac{\partial \mathcal{L}}{\partial z_w} = (\alpha + 1) z_w^{\alpha} - 2\lambda z_w = 0.
\end{equation}
When $z_w \neq 0$, it follows from \eqref{e26} that
\[
z_w^{\alpha - 1} = \frac{2\lambda}{\alpha + 1}.
\]
This implies that all non-zero $z_w$ are equal, denoted as $x$. Suppose that there are $k$ non-zero variables, then $k x^2 = 1$, which leads to $x = \frac{1}{\sqrt{k}}$. Substituting it into $V_2$ gives
\[
V_2 = k \Big( \frac{1}{\sqrt{k}} \Big)^{\alpha+1} = k^{1 - \frac{\alpha+1}{2}} = k^{\frac{1-\alpha}{2}}.
\]
For $0 < \alpha < 1$, $\frac{1-\alpha}{2}$ is positive, which implies that $V_2$ increases as $k$ increases. When $k = 1$, we have $V_2 = 1^{\frac{1-\alpha}{2}} = 1$. In this case, only one $z_w$ is nonzero and equals 1, while the other $z_w$ are zero. Thus, the minimum value of $V_2(Z)$ under the constraint $V_1(Z) = 1$ is 1 and \eqref{4inequality} reduces to \( \left( V_1(r) \right)^{\frac{\alpha + 1}{2}} \leq V_2(r) \). Hence, the inequality in the conclusion is established.
\begin{thm}\label{Theorem3}
Under the action of the controller \eqref{e8}, the system \eqref{b1} with the Hamiltonians \eqref{b2} and \eqref{b3} is globally finite-time stable, that is, the system will be stabilized in the equivalence class of the target
state $\left| {{\psi }_{f}} \right\rangle =\left| n-1 \right\rangle $ within a finite time.
\end{thm}
\textbf{Proof.}
Define \( g(t) = | \cos \phi_{n(n-1)}|^{\alpha_{n-1} + 1} \) as a function of phases. The asymptotic stability established in Theorem~\ref{Theorem2} ensures the existence of a time \( T_1 \) such that the amplitude of the target state satisfies \( r_n \geq \beta \), where \( \beta \in (0,1) \).
For \( t \ge T_1 \), the first-order time derivative \eqref{e9} of the Lyapunov function \eqref{e6} can be expressed as
\begin{equation}\label{V_dot_g1}
\dot{V} \leq -2 \beta k_{n-1} r_{n-1}^{\alpha_{n-1} + 1} g(t).
\end{equation}
When $n=2$, we let ${{K}_{f1}}=2\beta {{k}_{1}}$ and ${{\alpha }_{f1}}=\frac{{{\alpha }_{1}}+1}{2}$, and have
\begin{equation}\label{V_dot_finite1}
\dot{V} \leq -K_{f_1} (r_1^2)^{\alpha_{f_1}} g(t) = -K_{f_1} (V)^{\alpha_{f_1}} g(t).
\end{equation}
According to Theorem 4 in \cite{KUANG2021109327}, the system \eqref{b1} is finite-time stable. 

When \( n > 2 \), if the initial state does not satisfy \( r_{n-1} \geq r_j > 0 \) for all \( 1 \leq j \leq n-2 \), then the asymptotic stability established in Theorem~\ref{Theorem2} ensures that the inequality \( r_{n-1} \geq r_j > 0 \) holds after the time \( T_2 \). If the initial state satisfies \( r_{n-1} \geq r_j > 0 \) for all \( 1 \leq j \leq n-2 \), then \( T_2=T_0 \). Thus, for $t\ge \max ({{T}_{1}},{{T}_{2}})$, we have
\[
\begin{aligned}
\dot{V} &\leq -2 \beta k_{n-1} r_{n-1}^{\alpha_{n-1} + 1} g(t) \leq -\frac{2 \beta k_{n-1}}{n-1} g(t) \sum_{i=1}^{n-1} r_i^{\alpha_{n-1} + 1}.
\end{aligned}
\]
Applying Lemma~\ref{Lemma1} to the above equation, we obtain 
\[\dot{V} \leq -\frac{2 \beta k_{n-1}}{n-1} g(t) \Big( \sum_{j=1}^{n-1} r_j^2 \Big)^{\frac{\alpha_{n-1} + 1}{2}}.\]
Letting ${{K}_{f}}=\frac{2\beta {{k}_{n-1}}}{n-1}$ and ${{\alpha }_{f}}=\frac{{{\alpha }_{n-1}}+1}{2}$, we have
\begin{equation}\label{V_dot_finite3}
\dot{V} \leq  -K_{f_2} (V)^{\alpha_{f_2}} g(t).
\end{equation}
Since the relative phase does not affect the stability of the system, from \eqref{V_dot_finite1}, \eqref{V_dot_finite3}, and Theorem 4 in \cite{KUANG2021109327}, it follows that the system \eqref{b1} is finite-time stable for $t\ge \max ({{T}_{1}},{{T}_{2}})$. That is, the system \eqref{b1} is globally finite-time stable and stabilizes in the equivalence class of the target state \( \left| \psi_f \right\rangle = \left| n-1 \right\rangle \). 

Moreover, if the initial time satisfies \( t \geq \max (T_1, T_2) \), then the time required for exact convergence to the target state can be given by
\[
T_f \leq \frac{1}{K_{f_2} (1 - \alpha_{f_2})} V(0)^{1 - \alpha_{f_2}}.
\]

\section{Numerical Simulations}
\begin{figure}[htbp]
\centering
\includegraphics[width=0.8\linewidth]{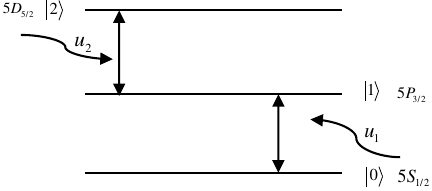}  
\caption{Configuration of a there-level ladder rubidium atom, where \( 5S_{1/2} \), \( 5P_{3/2} \), and \( 5D_{5/2} \) correspond to \( |0\rangle \), \( |1\rangle \), and \( |2\rangle \), respectively, and the control fields $u_1$ and $u_2$ are exerted to the system shown in the figure.}
\label{fig1}
\end{figure}

We consider a ladder three-level rubidium atomic system with the energy levels and the control fields shown in Fig.~\ref{fig1}, whose Hamiltonians can be expressed as
$H_0 =\left[ \begin{smallmatrix}
0 & 0 & 0 \\
0 & 1 & 0 \\
0 & 0 & 2
\end{smallmatrix}\right]$, 
$H_1 = \left[ \begin{smallmatrix}
0 & -i & 0\\
i & 0 & 0\\
0 & 0 & 0
\end{smallmatrix}\right]$, and
$H_2 = \left[ \begin{smallmatrix}
0 & 0 & 0\\
0 & 0 & -i\\
0 & i & 0
\end{smallmatrix}\right]$.

\begin{figure}[tbp]
\centering
\includegraphics[width=0.8\linewidth]{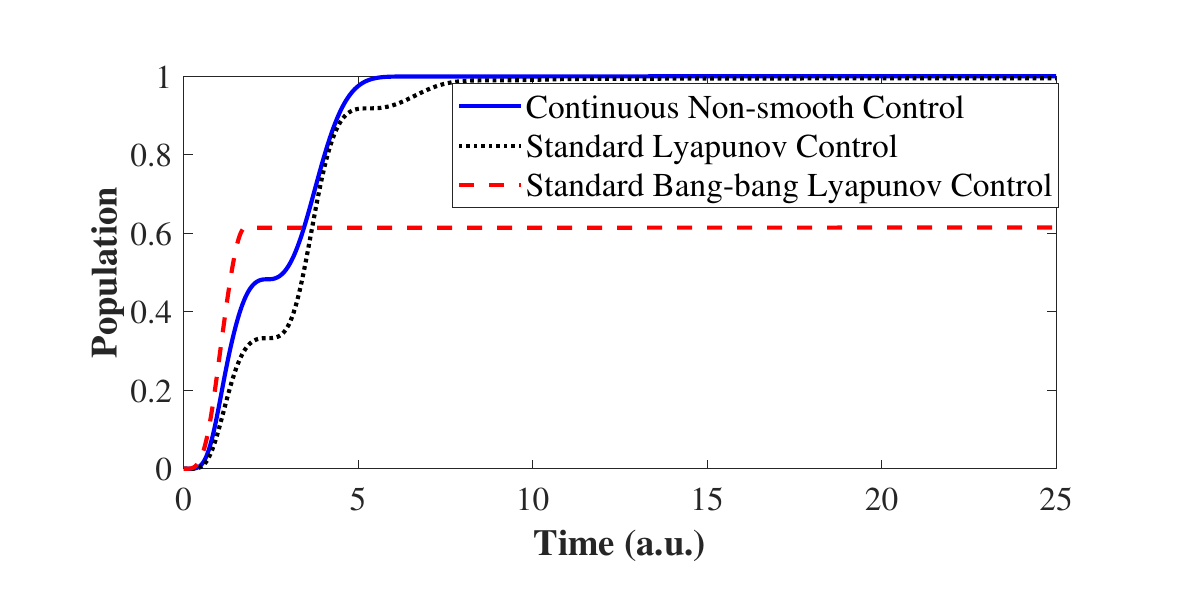}  
\caption{For the initial state \( |\psi(0)\rangle = |0\rangle \) and the target state \( |\psi_f\rangle = |2\rangle \), the population evolution curves of the ladder quantum system under the continuous non-smooth controls, the standard Lyapunov controls, and the standard bang-bang Lyapunov controls.}
\label{fig3}
\end{figure}

\begin{figure}[tbp]
\centering
\includegraphics[width=0.8\linewidth]{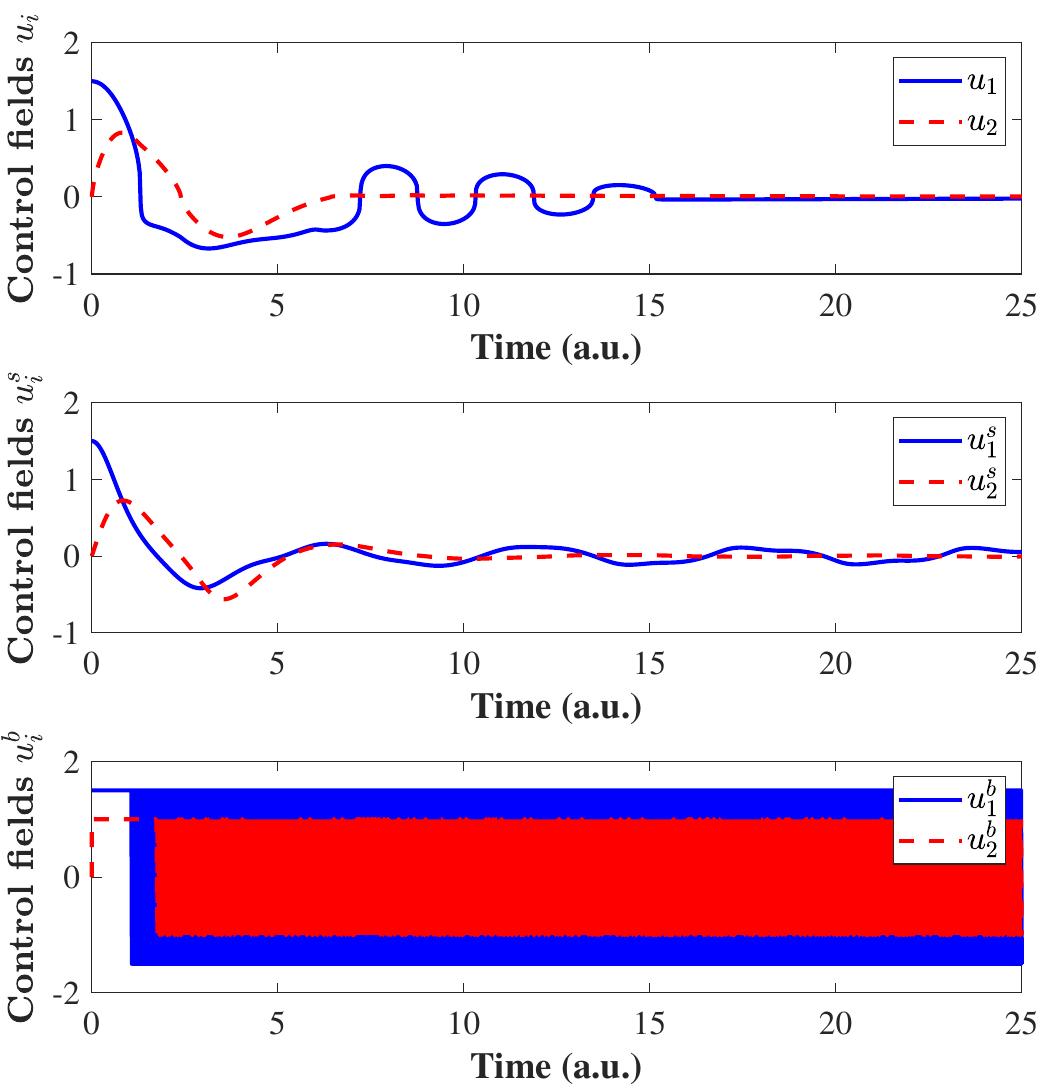}  
\caption{For $|\psi(0)\rangle = |0\rangle$, the continuous non-smooth controls $u_1$ and $u_2$, the standard Lyapunov controls $u_1^s$ and $u_2^s$, and the standard bang–bang Lyapunov controls $u_1^b$ and $u_2^b$.}
\label{fig4}  
\end{figure}

We assume that the initial state is \( |0\rangle = [1, 0, 0]^T \) and the target state is \( |2\rangle = [0, 0, 1]^T \). Generally, the system can naturally enter the finite-time effective region when the condition \( t \geq \max (T_1, T_2) \) in Theorem~\ref{Theorem3} is satisfied. However, to reduce the arrival time, we can increase \( u_1 \) such that the condition \( r_2 \geq r_1 > 0 \) is satisfied earlier and accordingly the finite-time effective stability region is entered earlier. Actually, this can be achieved by adjusting the parameters in the control law~\eqref{e8}, e.g., setting \( k_1 \) larger than \( k_2 \) and setting \( \alpha_1 \) smaller than \( \alpha_2 \). In the simulation, we choose the controller parameters as \( k_1 = 1.5 \), \( k_2 = 1 \), \( \alpha_1 = 1/3 \), \( \alpha_2 = 2/3 \), and have
\begin{equation}\label{u1_u2}
\begin{cases}
u_1 = 1.5 \operatorname{sign}(r_1 \cos \phi_{21}) \left| r_1 \cos \phi_{21} \right|^{1/3} \\
u_2 = \operatorname{sign}(r_2 \cos \phi_{32}) \left| r_2 \cos \phi_{32} \right|^{2/3}.
\end{cases}
\end{equation}
We also conduct simulation experiments under the standard Lyapunov control \( u^{s} \) \cite{CONG200728} and the standard bang–bang Lyapunov control \( u^{b} \) \cite{KUANG2017164} to compare their performance, where the standard Lyapunov controls are given by \( u_1^s = 1.5 r_1 \cos \phi_{21} \) and \( u_2^s = r_2 \cos \phi_{32} \), and the standard bang–bang Lyapunov controls are \( u_1^b = 1.5 \operatorname{sign}(r_1 \cos \phi_{21}) \) and \( u_2^b = \operatorname{sign}(r_2 \cos \phi_{32}) \). The simulation results are shown in Fig. \ref{fig3} and Fig. \ref{fig4}.

According to Fig. \ref{fig3} and the simulation data, the system stabilizes in the target state $| 2 \rangle$ after about $t_f\approx 16.6814$ a.u. However, at this moment, the populations of the target state under the standard Lyapunov controls $u_j^{s}$ $(j=1,2)$ and the standard bang–bang Lyapunov controls $u_j^{b}$ $(j=1,2)$ are 0.9944 and 0.6144, respectively. As shown in Fig. \ref{fig4}, the proposed control law in this paper is indeed continuous and non-smooth, whereas the standard Lyapunov control law is smooth, and the standard bang-bang Lyapunov control law is discontinuous with chattering. The simulation results demonstrate that under the action of the controller \eqref{u1_u2}, the system \eqref{b1} with the Hamiltonians \eqref{b2} and \eqref{b3} stabilizes in the equivalence class of the target state $| {{\psi }_{f}} \rangle =| 2 \rangle $ within a finite time.

\begin{figure}[htbp]
\centering
\includegraphics[width=0.8\linewidth]{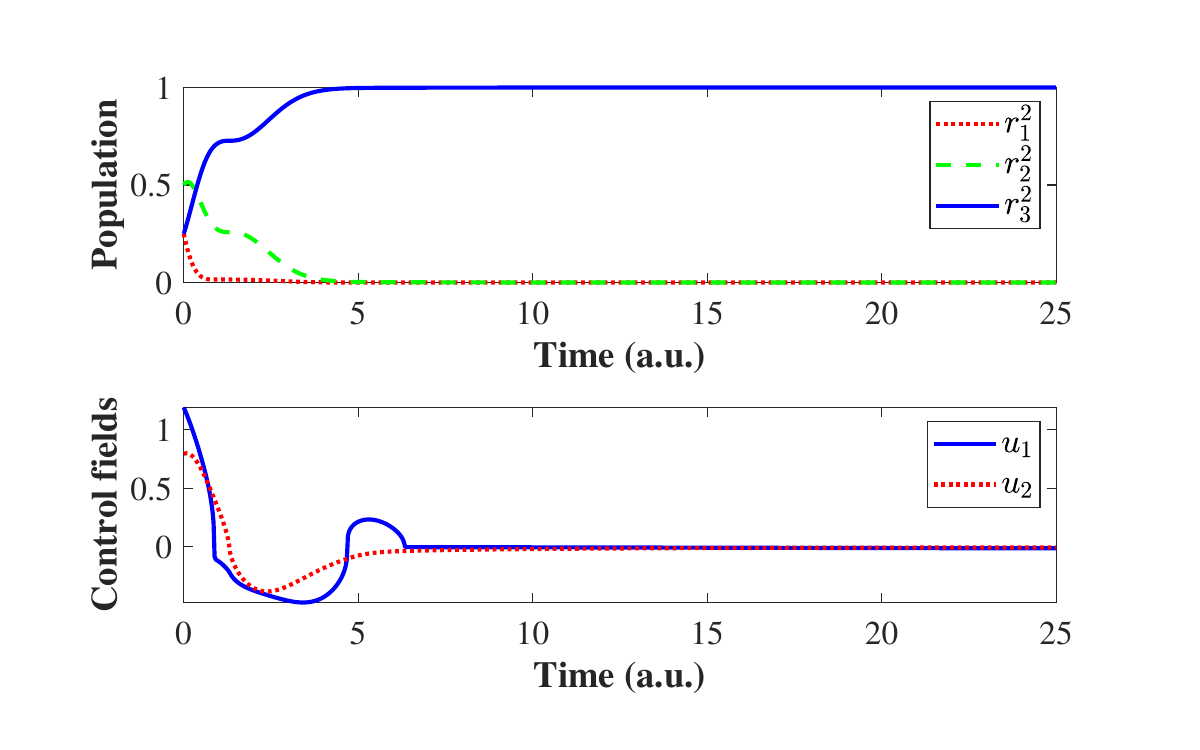}    
\caption{For the initial state \( | \psi (0) \rangle = [ \frac{1}{2}, \frac{\sqrt{2}}{2}, \frac{1}{2} ]^T \) and the target state \( | \psi_f \rangle = | 2 \rangle \), the population evolution curves and the continuous non-smooth controls \( u_1 \) and \( u_2 \) in the ladder quantum system.}
\label{fig5}                                 
\end{figure}

Next, we select the initial state as $| \psi (0) \rangle ={{[\frac{1}{2},\frac{\sqrt{2}}{2},\frac{1}{2}]}^{T}}$ and the target state as $| 2 \rangle ={{[0,0,1]}^{T}}$. The simulation results are shown in Fig.~\ref{fig5}. It is observed that the condition ${{r}_{2}}\ge {{r}_{1}}>0$ always holds throughout the whole evolution process. Let $\beta =\frac{1}{2}$ in \eqref{V_dot_g1}. Fig.~\ref{fig5} and the simulation data show that the stabilization time of the system is $t_f\approx 8.4940$ a.u. According to Theorem 4 in \cite{KUANG2021109327}, the stabilization time corresponding to the initial state $| \psi (0) \rangle ={{[\frac{1}{2},\frac{\sqrt{2}}{2},\frac{1}{2}]}^{T}}$ satisfies ${{t}_{f}}\approx 8.4940\text{ a}\text{.u}\text{.}<\frac{6}{{{K}_{2}}(1-{{\alpha }_{2}})}V{{(| \psi (0) \rangle )}^{\frac{1-{{\alpha }_{2}}}{2}}}=17.1573\text{ a}\text{.u}\text{.}$, which is consistent with the theoretical result.

\section{Conclusion}
We proposed a novel continuous non-smooth control law and achieved the finite-time convergence of any $n$-level ladder quantum system to an eigenstate of its internal Hamiltonian under the action of $n-1$ external control fields. Despite the presence of numerous non-Lipschitz regions, we showed that the high-dimensional quantum control system has a unique solution via the Filippov's solution and the LaSalle's invariance principle. The asymptotic and finite-time stability of the $n$-level control system was further validated by using the Lyapunov stability and the finite-time Lyapunov stability theories. Numerical simulations on a three-level ladder rubidium atom verify the effectiveness of the proposed control law. Future work will focus on optimizing the controller parameters $k$ and $\alpha$ to improve performance, as well as extending finite-time convergence control to $n$-level quantum systems with alternative structures.

\section*{Appendix}\label{Appendix}
We list several necessary lemmas used in this paper here.

\begin{lemma}\label{Lemma2} \cite{han2016theory} Consider the Cauchy problem of differential equation
$\dot{x} = f(x), \quad x_0 \in \mathbb{R}^n$,
where \( f(x) \) is bounded, but may be discontinuous. Its Filippov extension is
$F(x) = \bigcap_{\delta > 0} \text{clco} \, f(x + \delta B)$.
Then, the Cauchy problem has a differential inclusion solution 
$\dot{x} \in F(x), x_0 \in \mathbb{R}^n$. Moreover, let \( x(t, x_0) \) be any solution of the Filippov extension, then $\dot{x}(t, x_0) = f(x(t, x_0))$ if \( f(x) \) is continuous at \( x(t, x_0) \).
\end{lemma}

\begin{lemma}\label{Lemma3}
\cite{dixon2001nonlinear} Suppose \( f(t): \mathbb{R}^+ \to \mathbb{R} \) is a differentiable function and that \(\lim_{t \to \infty} f(t) = C\) exists. If the derivative of the function $f(t)$ can be expressed as the sum of two functions, i.e.,
$f'(t) = g_1(t) + g_2(t)$, where \( g_1(t) \) is uniformly continuous and \( \lim_{t \to \infty} g_2(t) = 0 \), then
$\lim_{t \to \infty} f(t) = 0$ and $\lim_{t \to \infty} g_1(t) = 0$. 
\end{lemma}

\begin{lemma}\label{Lemma4} \cite{Bhat2005} Let $V_1$ and $V_2$ be continuous real-valued functions defined on $\mathbb{R}^{n-1}$ and $V_1$ be positive definite. Suppose that $V_1$ and $V_2$ are homogeneous of degrees $l_1 > 0$ and $l_2 > 0$ with respect to $\delta_\epsilon^d$ respectively. Then, for every $r \in \mathbb{R}^{n-1}$, the following holds
\[
\begin{aligned}
&\Big(\min_{\{z:V_1(z)=1\}}V_2(z)\Big)\big(V_1(r)\big)^{\frac{l_2}{l_1}}
\leq V_2(r)\\
\leq&\;\Big(\max_{\{z:V_1(z)=1\}}V_2(z)\Big)\big(V_1(r)\big)^{\frac{l_2}{l_1}}.
\end{aligned}
\]
\end{lemma}




\end{document}